\documentclass[11pt]{article}
\usepackage{amssymb}
\usepackage{amsmath}
\usepackage{amsfonts}
\usepackage{graphicx}
\usepackage{epstopdf}
\newtheorem{theorem}{Theorem}

\newtheorem{definition}{Definition}

\newtheorem{lemma}{Lemma}

\newtheorem{remark}{Remark}

\numberwithin{equation}{section}

\begin{document}
\title{Non-sharpness of the Morton-Franks-Williams inequality}
\author{{Keiko Kawamuro}\thanks{The author acknowledges partial support
from NSF grants DMS-0405586 and DMS-0306062}}
\date{January 28 2006.}
\maketitle
\begin{abstract}
We give (Theorem \ref{thmA}) conditions on a knot on which the
Morton-Franks-Williams inequality is not sharp. As applications,
we show infinitely many examples of knots where the
inequality is not sharp and also prove (by giving examples) that the
deficit of the inequality can be arbitrarily large.
\end{abstract}

\section{Introduction.}\label{MFW}
The Morton-Franks-Williams (MFW) inequality \cite{Morton},
\cite{FW}, is one of the few tools available in knot theory to
estimate the minimal braid index of a knot or a link. %

\begin{figure}[htpb!]
\begin{center}\includegraphics [height=20mm]{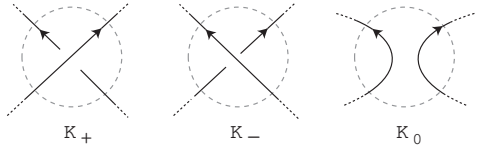}
\end{center}
\caption{Local views of $K_+, K_-, K_0$.}\label{+-0}
\end{figure}

To state the MFW inequality, let $K$ be an oriented knot or link
projected on a plane. Focus on one crossing of $K$ with sign
$\varepsilon$. Denote $K_\varepsilon := K$ and let
$K_{-\varepsilon}$ (resp. $K_0$) be the closed braid obtained from
$K_\varepsilon$ by changing the the crossing to the
opposite sign $-\varepsilon$ (resp. resolving the crossing),
see Figure \ref{+-0}. %

The \textit{HOMFLYPT polynomial} $P_{K}(v,z)$ of $K$ satisfies the
following relations (for any choice of a crossing):
\begin{eqnarray}
\frac{1}{v} P_{K_{+}} - v P_{K_{-}} &=& z P_{K_{0}}.  \label{skein1} \\
P_{\mathrm{unknot}} &=&1.  \notag
\end{eqnarray}
\textbf{The Morton-Franks-Williams inequality.} \ \ {\it Let
$d_{+}$ and $d_{-}$ be the maximal and minimal degrees of the
variable $v$ of $P_{K}(v,z)$. If a knot type ${\cal K}$ has a
closed braid representative $K$ with braid index $b_K$ and
algebraic crossing number $c_K$, then we have
$$
c_K - b_K +1 \leq d_{-}\leq d_{+}\leq c_K + b_K -1.
$$
As a corollary,
\begin{equation}
\frac{1}{2}(d_{+}-d_{-})+1 \leq b_K,  \label{lower bound}
\end{equation}%
giving a lower bound for the {\em braid index} $b_{\cal K}$ of
${\cal K}$.}

This inequality was the first known result of a general nature
relating to the computation of braid index, and it appeared to be
quite effective.  Jones notes, in \cite{Jones-1}, that on all but
five knots ($9_{42}, 9_{49}, 10_{132}, 10_{150}, 10_{156}$) in the
standard knot table, up to crossing number $10$, the MFW
inequality is sharp.
Furthermore it has been known that the inequality is sharp on all
torus links, closed positive $n$-braids with a full twist \cite{FW}, $2$-bridge
links \cite{Murasugi}, fibred  and alternating links \cite{Murasugi}.

However, the MFW inequality is not as strong as it appears to be
as above. In fact, in Theorem \ref{deficit-thm} and Theorem
\ref{BM-thm}, we give infinitely many examples of prime knots and
links on which the MFW inequality is not sharp and is arbitrarily
far away from being an equality. All these examples are obtained
as corollaries of our main result Theorem \ref{thmA}, in which we
give one reason to explain non-sharpness of the MFW inequality.
The main idea is to find knots $K_\alpha$ of known braid index $=
b$ which have a distinguished crossing such that, after changing
that crossing to each of the other two possibilities in Figure
\ref{+-0}, giving knots or links $K_\beta$ and $K_\gamma$, it is
revealed that $K_\beta$ and $K_\gamma$ each has braid index $< b.$
Thanks to this theorem one can observe ``accumulation'' of
deficits by looking at the distinguished crossings which
contribute to deficits (for detail, see the proof of Theorem
\ref{deficit-thm}).

\textbf{Acknowledgment. }
This paper is part of the author's
ongoing work toward her Ph.D. thesis. She is grateful to her
advisor, Professor Joan Birman, for very much thoughtful advice
and for her encouragement.
She also wishes to thank Professor
William Menasco, who told her about the Birman-Menasco diagram,
discussed in Section \ref{section-BM} and the associated
conjecture, when she visited him at SUNY Buffalo in July 2004, to
Professor Walter Neumann for helpful suggestions and to Professor 
Alexander Stoimenow for sending a preprint. Finally, she especially
thanks Professor Mikami Hirasawa, who shared many creative ideas and
results about fibred knots including the definition and properties of
the enhanced Milnor number.

\section{One reason for non-sharpness of the MFW inequality.}\label{reason}
In this section we give sufficient conditions (Theorem \ref{thmA})
for a closed braid on which the MFW inequality is not sharp.
Then we exhibit examples of prime links on which the deficit
of the inequality can be arbitrary large.

Let $b_{\cal K}$ be the braid index of knot type ${\cal K}$,
that is the smallest integer $b_{\cal K}$ such that ${\cal K}$ can
be represented by a closed $b_{\cal K}$-braid. Let $b_K, c_K$ denote the
braid index and the algebraic crossing number of a braid
representative $K$ of ${\cal K}$.

\begin{definition}Let
\begin{equation*}
D_{\cal K} : = b_{\cal K} - \frac{1}{2}(d_{+}-d_{-})-1
\end{equation*}%
be the difference of the numbers in $(\ref{lower bound}),$ i.e.,
of the actual braid index and the lower bound for braid index.
Call $D_{\cal K}$ the {\em deficit} of the MFW inequality for
${\cal K}$.
\end{definition}
If $D_{\cal K} =0,$ the MFW inequality is sharp on ${\cal K}$. If
$K$ is a braid representative of ${\cal K}$ let $D_K^+ \ := \ (c_K
+ b_{K} -1)-d_{+}$ and $D_K^- \ := \ d_{-}-(c_K - b_{K} +1).$ When
$b_K = b_{\cal K},$ we have
\begin{equation}\label{D+D}
 D_{\cal K} =\frac{1}{2}(D_K^{+} + D_K^{-}).
\end{equation}%
Note that $D_K^{\pm}$ depends on the choice of braid
representative $K$, but the deficit $D_{\cal K}$ is independent
from the choice.

\begin{theorem}\label{thmA}
Assume that $K$ is a closed braid representative of ${\cal K}$ with
$b_K = b_{\cal K}$. Focus on one crossing of $K$ and construct
$K_+, K_-, K_0$ $($one of the three must be $K)$. Let $\alpha,
\beta, \gamma \in \{ +, -, 0 \}$ and assume that $\alpha, \beta,
\gamma$ are mutually distinct. If $K_\alpha = K$ and if positive
$($resp. negative$)$ destabilization is applicable $p$-times
$($resp. $n$-times$)$
to each of $ K_\beta$ and $K_\gamma$, then
\begin{eqnarray}
D_K^+ & \geq & 2p,    \label{less}     \\
( \text{resp.} \quad D_K^- &\geq & 2n. ) \label{more}%
\end{eqnarray}%
i.e., by {\em (\ref{D+D})} the MFW inequality is not sharp on ${\cal K}$
if $p+n>0$.
\end{theorem}

Here is a lemma to prove Theorem \ref{thmA}.
\begin{lemma}\label{lemma-for-thmA}
Let $K$ be a closed braid. Choose one crossing, and construct
$K_+, K_-, K_0$ $($one of the three must be $K)$. We have
\begin{eqnarray}
d_{+}(P_{K_{+}}) &\leq & \max \{d_{+}(P_{K_{-}})+2, \quad d_{+}(P_{K_{0}})+1\} \label{dplus+} \\
d_{+}(P_{K_{-}}) &\leq & \max \{d_{+}(P_{K_{+}})-2, \quad d_{+}(P_{K_{0}})-1\} \label{dplus-} \\
d_{+}(P_{K_{0}}) &\leq & \max \{d_{+}(P_{K_{+}})-1, \quad d_{+}(P_{K_{-}})+1\} \label{dplus0}
\end{eqnarray}%
and%
\begin{eqnarray*}
d_{-}(P_{K_{+}}) &\geq & \min \{d_{-}(P_{K_{-}})+2, \quad d_{-}(P_{K_{0}})+1\} \\
d_{-}(P_{K_{-}}) &\geq & \min \{d_{-}(P_{K_{+}})-2, \quad d_{-}(P_{K_{0}})-1\} \\
d_{-}(P_{K_{0}}) &\geq & \min \{d_{-}(P_{K_{+}})-1, \quad d_{-}(P_{K_{-}})+1\}. \end{eqnarray*}
\end{lemma}

\textbf{Proof of Lemma \ref{lemma-for-thmA}. } By
(\ref{skein1}), we have $P_{K_{+}}=v^{2}P_{K_{-}} + vz P_{K_{0}}.$
Thus, $d_+(P_{K_+})=d_+(v^2 P_{K_-} + vz P_{K_0}) \leq \max \{ d_+
( v^2 P_{K_-}), \ d_+ ( vz P_{K_0})\}$ and we obtain
(\ref{dplus+}). The other results follow similarly.

 \hfill $\Box$ %

Table (\ref{table1}) shows the changes of $c_K$, $b_K$, $c_K - b_K + 1$ and
$c_K + b_K - 1$
under stabilization and destabilization of a closed braid.
\begin{equation}
\begin{tabular}{|l|l|l|l|l|}
\hline & $c_K$ & $b_K$ & $c_K-b_K+1$ & $c_K+b_K-1$ \\ \hline $+$ stabilization
& $+1$ & $+1$ & \multicolumn{1}{|c|}{$0$} &
\multicolumn{1}{|c|}{$+2$} \\ \hline $+$ destabilization & $-1$ &
$-1$ & \multicolumn{1}{|c|}{$0$} & \multicolumn{1}{|c|}{$-2$} \\
\hline $-$ stabilization & $-1$ & $+1$ &
\multicolumn{1}{|c|}{$-2$} & \multicolumn{1}{|c|}{$0$} \\ \hline
$-$ destabilization & $+1$ & $-1$ & \multicolumn{1}{|c|}{$+2$} &
\multicolumn{1}{|c|}{$0$} \\ \hline
\end{tabular}
\label{table1}
\end{equation}
Note that $c_K$ and $b_K$ are invariant under braid isotopy and exchange moves.

\textbf{Proof of Theorem \ref{thmA}.  } Suppose that $K = K_\alpha
= K_+.$  Suppose we can apply positive destabilization $k$-times
($k\geq p$) to $K_{-}$. Let $\tilde{K}_{-}$ denote the closed
braid by the destabilization. Then we have:%
\begin{eqnarray}
d_{+}(P_{K_-}) + 2 &=& d_{+}(P_{\tilde{K}_-}) + 2  \notag \\
&\leq &( c_{\tilde{K}_-} + b_{\tilde{K}_-} -1 ) + 2  \notag \\
&=& \{(c_{K_-} + b_{K_-} -1) -2k \} + 2  \label{5eq} \\
&=&   (c_{K_+} -2) + b_{K_+} -1 -2k + 2   \notag \\
&=&   (c_{K_+} + b_{K_+} -1) -2k =  (c_K + b_K -1) -2k. \notag
\end{eqnarray}
The first equality holds since $K_{-}$ and $\tilde{K}_-$ have the
same knot type. The first inequality is the MFW inequality. The
second equality follows from Table (\ref{table1}).

Similarly, if we can apply positive destabilization $l$-times $(l
\geq p)$ to $K_{0},$ and obtain $\tilde{K_0}$, we have
\begin{eqnarray}
d_{+}(P_{K_{0}})+1 &=& d_+ (P_{\tilde{K_0}}) + 1  \notag \\
&\leq& (c_{\tilde{K_0}} + b_{\tilde{K_0}} -1 ) + 1  \notag \\
&=&    (c_{K_0} + b_{K_0} -1 - 2l ) + 1  \label{5eq'} \\
&=&    (c_{K_+} -1 ) + b_{K_+} -1 - 2l + 1 \notag \\
&=&    (c_{K_+} + b_{K_+} -1) -2l =  (c_K + b_K -1) -2l. \notag
\end{eqnarray}

By (\ref{dplus+}), (\ref{5eq}) and (\ref{5eq'}) we get
\begin{eqnarray*}
d_+ (P_K) &=& d_+ (P_{K_+}) \leq \max \{d_{+}(P_{K_{-}})+2, \quad d_{+}(P_{K_{0}})+1\} \\
          &\leq& (c_K + b_{\cal K} -1) - \min \{ 2k, 2l \},
\end{eqnarray*}
i.e., $D_K^+ \geq \min \{ 2k, 2l \} \geq 2p.$
When $K_\alpha = K_-$ or $K_\alpha = K_0$, the same arguments work
(use (\ref{dplus-}) or (\ref{dplus0}) for these cases in the place
of (\ref{dplus+})) and we get (\ref{less}).

The other inequality (\ref{more}) also holds by the identical
argument.
\hfill $\Box$ %

\begin{figure}[htpb!]
\begin{center}
\includegraphics [height=55mm]{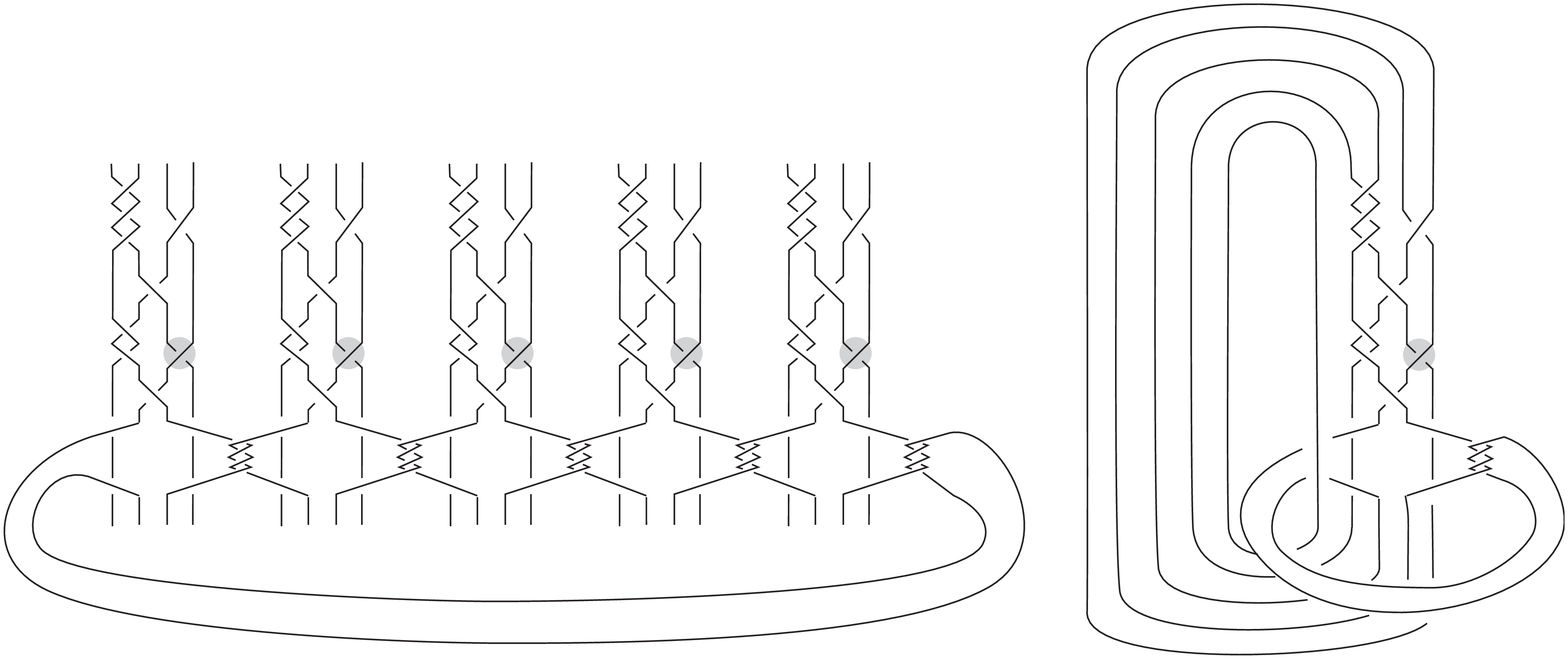}
\end{center}
\caption{Prime link $(9_{42})^5$ and $2$-component link.} \label{942}%
\end{figure}

\begin{theorem}\label{th942}
Knot type ${\cal K}=9_{42}$ has a braid representative $K=K_+$
$($see Figure \ref{942}$)$ satisfying the sufficient condition in
Theorem \ref{thmA}.
\end{theorem}
\begin{figure}[htpb!]
\begin{center}\includegraphics [height=40mm]{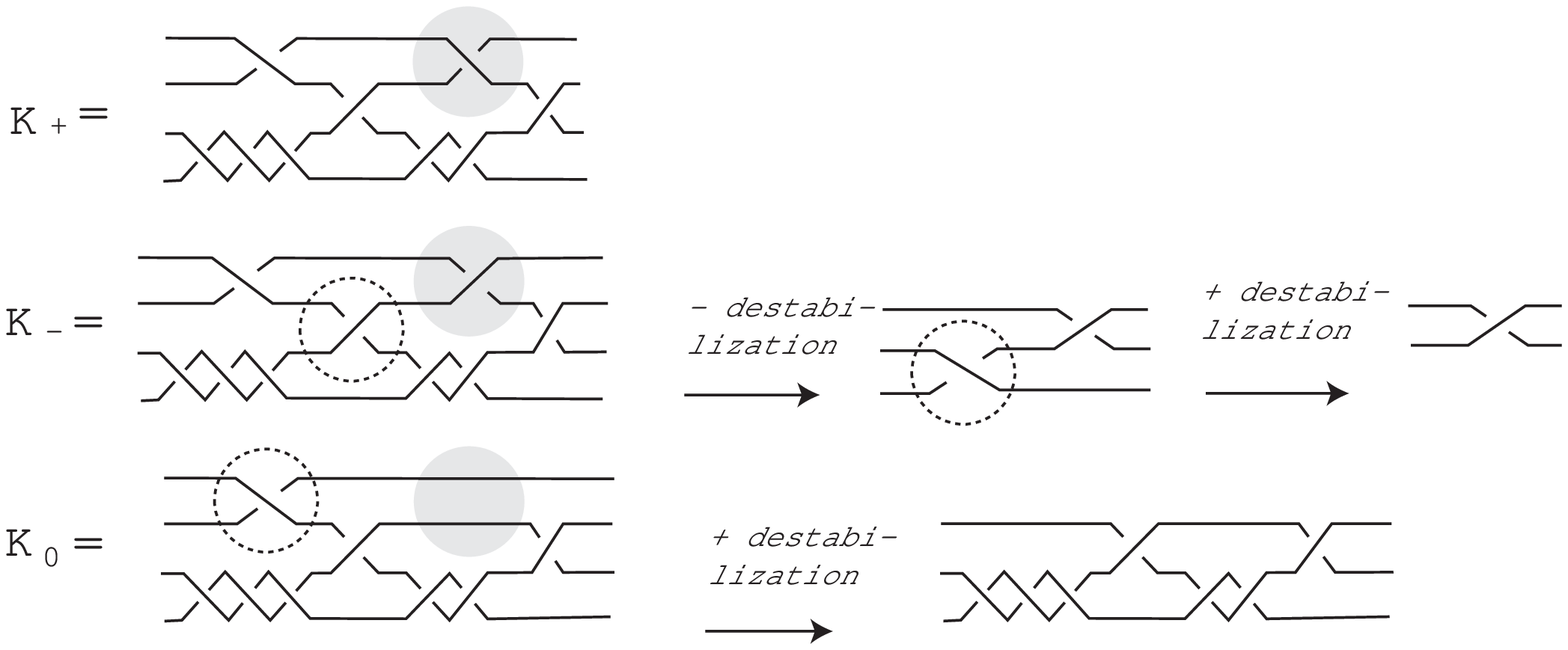}
\end{center}
\caption{Knot $9_{42}$ satisfies the conditions of Theorem
\ref{thmA}}\label{skein942}
\end{figure}
\textbf{Proof of Theorem \ref{th942}. }
It is known that $9_{42}$ has braid index $=4$ and deficit $D_{9_{42}}=1.$
Let $K = K_+$ be its braid representative of the minimal braid index
as in Figure \ref{skein942}.
Construct $K_-, K_0$ by changing the shaded crossing.
Sketches show that both $K_-, K_0$ can be positively
destabilized. Thus by Theorem \ref{thmA}, $D_{K}^+ \geq 2$ and
$D_{9_{42}}\geq 1.$
\hfill $\Box$ %

\begin{theorem}\label{deficit-thm}
For any positive integer $n,$ there exists a prime link whose
deficit is $\geq n.$%
\end{theorem}

\textbf{Proof of Theorem \ref{deficit-thm}. } We prove the theorem
by exhibiting examples. For $n \in \mathbb{N}$ let $(9_{42})^n$ be
the closure of $n$-copies of $9_{42}$ linked each other by two
full twists as in the left sketch of Figure \ref{942}.
Since the braid index $b_{9_{42}} = 4$ and $(9_{42})^n$ is
an $n$-component link, we know the braid index of $(9_{42})^n$ is
$4n.$ This construction gives a braid representative with
$4n$-strands and $n$ distinguished (shaded in the left sketch)
crossings.

In the following we will see that each of the shaded crossing
contributes to the deficit.

Let ${\cal K}:=(9_{42})^2$ and let $K$ be the braid representative
of ${\cal K}$ as in Figure \ref{942}. Let $K_{- -}, K_{- 0}, K_{0
-}, K_{0 0}$ be the links obtained from $K$ by changing the two
shaded crossings. We repeat the discussion of the proof of Theorem
\ref{thmA}: We have
\begin{eqnarray*}
d_+(P_{K_{- -}}) + (2+2) &=& d_{+}(P_{\tilde{K}_{- -}}) + 4   \\ %
&\leq &( c_{\tilde{K}_{- -}} + b_{\tilde{K}_{- -}} -1 ) + 4   \\
&=& \{(c_{K_{- -}} + b_{K_{- -}} -1) -2\cdot 2 \} + 4  \\
&=&   (c_{K} -4) + b_{K} -1 -2\cdot 2 + 4   \\
&=&   (c_{K} + b_{K} -1) -2\cdot 2.
\end{eqnarray*}
Similarly,
\begin{eqnarray*}
 d_+(P_{K_{- 0}})+(2+1) &\leq & (c_{K} + b_{K} -1) -2\cdot 2, \\
 d_+(P_{K_{0 -}})+(1+2) &\leq & (c_{K} + b_{K} -1) -2\cdot 2, \\
 d_+(P_{K_{0 0}})+(1+1) &\leq & (c_{K} + b_{K} -1) -2\cdot 2.
\end{eqnarray*}
Thus,
\begin{eqnarray*}
d_+ (P_K) &=& \max \{d_+(P_{K_{- -}})+4, \  d_+(P_{K_{- 0}})+3, \ %
                     d_+(P_{K_{0 -}})+3, \  d_+(P_{K_{0 0}})+2 \} \\ %
          &\leq & (c_K + b_{\cal K} -1) - 2\cdot 2
\end{eqnarray*}
and $D_{\cal K} \geq \frac{1}{2} D_K^+ \geq \frac{1}{2}(2\cdot 2) = 2.$

Similar arguments work when ${\cal K}=(9_{42})^n$ for $n\geq 3$
and we have
 $D_{(9_{42})^n} \geq \frac{1}{2} D_{(9_{42})^n}^+ \geq \frac{1}{2} (2\cdot n)
 \geq n.$

Since the $2$-component link of the right sketch is hyperbolic
\cite{N}, by \cite{NZ} we can conclude that $(9_{42})^n$'s are all prime
except for finitely many cases.

\hfill $\Box $%

\begin{remark}
{\rm By taking the connected sum of knots on which the MFW
inequality is non-sharp, one can also construct examples of
(non-prime) knots with arbitrarily large deficits. This fact
follows not only by Theorem \ref{thmA} but also by
the definition of HOMFLYPT polynomial (\ref{skein1})
and the additivity of braid indices under connected sums
\cite{BM4}. }
\end{remark}

\section{The Birman-Menasco block and strand diagram.}\label{section-BM}
\begin{figure}[htpb!]
\begin{center}
\includegraphics [height=30mm]{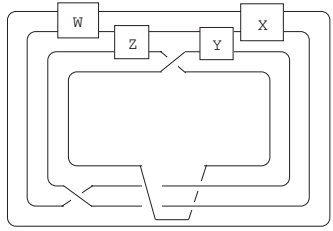}
\end{center}
\caption{The Birman-Menasco diagram $BM_{x,y,z,w}$.} \label{menasco}
\end{figure}

In this section as an application of Theorem \ref{thmA} we study
another infinite class of knots including all the Jones' five
knots on which the MFW inequality is not sharp.
We call the block-strand diagram (see \cite{MTWS-I} for definition)
of Figure \ref{menasco} the
Birman-Menasco (BM) block-strand diagram.
\begin{definition}
Let $BM_{x,y,z,w}$, where $x,y,z,w \in \mathbb{Z}$, be the knot
$($or the link$)$ type which is obtained by assigning $x$-half
positive twists $($resp. $y, z, w )$ to the braid block $X$
$($resp. $Y, Z, W)$ of the BM diagram.
\end{definition}

Recall that on all but only five knots ($9_{42}, 9_{49}, 10_{132},
10_{150}, 10_{156}$) up to crossing number $10$ the MFW inequality
is sharp. An interesting property of the BM diagram is that it
carries all the five knots. Namely, we have
$9_{42}=BM_{-1,1,-2,-1}=BM_{-1,-2,-2,2}, \ $
$9_{49}=BM_{-1,1,1,2}\ $, $10_{132} = BM_{-1,-2,-2,-2}, \ $
$10_{150} = BM_{3,-2,-2,2} = BM_{-1,2,-2,2} = BM_{-1,-2,2,2} =
BM_{-1,1,2,-1} =BM_{3,1,-2,-1},$ and $10_{156}=BM_{-1,1,1,-2}$.

We have the following theorem, which was conjectured informally by Birman and Menasco: %
\begin{theorem}\label{BM-thm}
There are infinitely many $(x,y,z,w)$'s such that
the MFW inequality is not sharp on $BM_{x,y,z,w}$.
\end{theorem}
We need lemmas to prove Theorem \ref{BM-thm}.
\begin{lemma}\label{D+}
We have $D_{BM_{x,y,z,w}}^+ \geq 2$.
\end{lemma}

\begin{figure}[htp]
\begin{center}
\includegraphics [height=120mm]{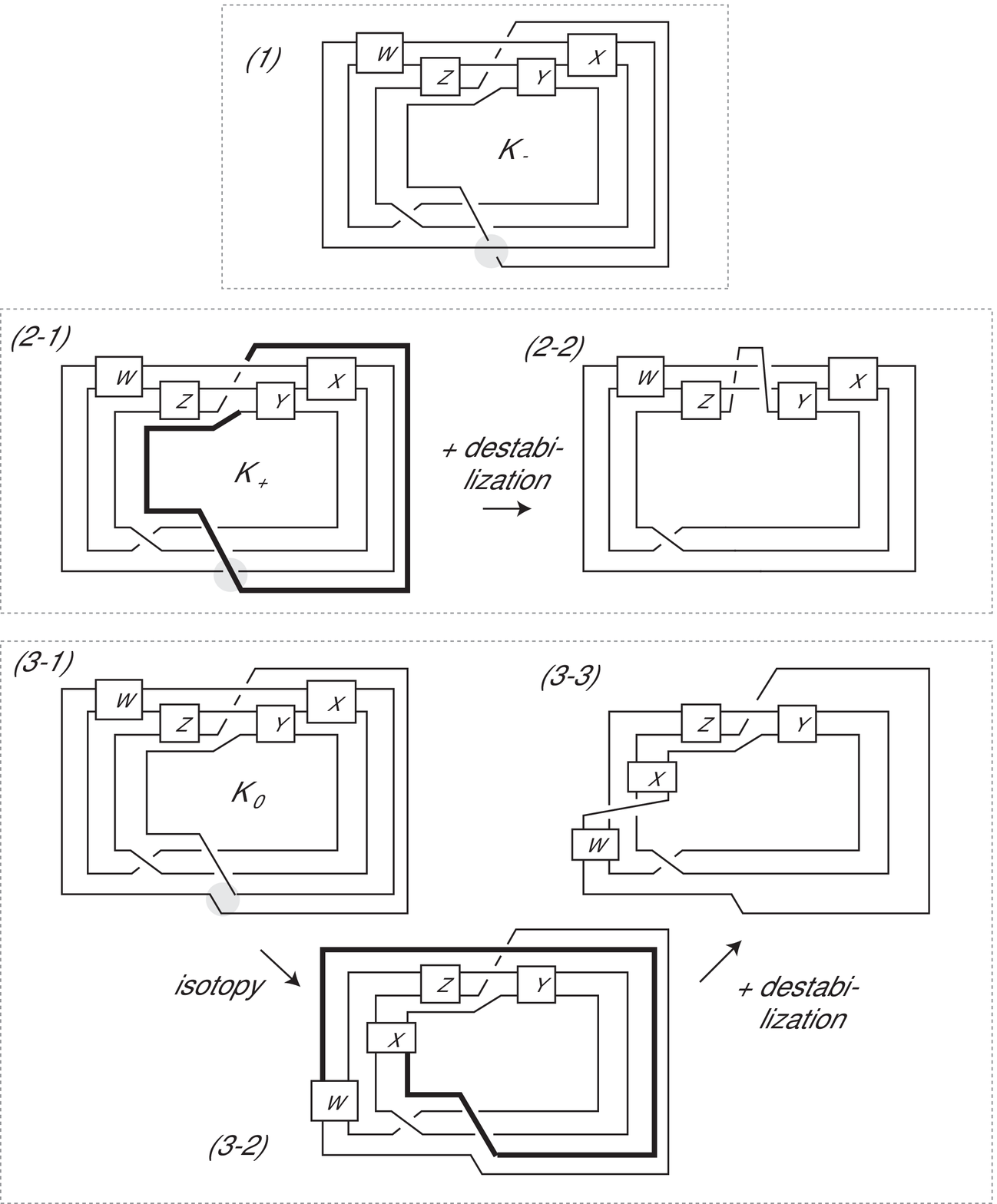}
\end{center}
\caption{} \label{menasco2}
\end{figure}

\textbf{Proof of Lemma \ref{D+}.} \ \
Change the BM diagram into the diagram in sketch (1) of Figure \ref{menasco2}%
by braid isotopy and denote it by $K$. Focus on the crossing
shaded in the sketch (1).  Regard $K = K_-.$ We can apply positive
destabilization once to $K_+$ and obtain the diagram in sketch
(2-2). We also can apply positive destabilization once to $K_0$ as
we can see in the passage sketch $($3-1$)$ $\Rightarrow$
$($3-2$)$ $\Rightarrow$ $($3-3$)$. Therefore by Theorem \ref{thmA}
we have $D_{BM_{x,y,z,w}}^+ \geq 2$ for any $(x,y,z,w)$.

\hfill $\Box$ %

It remains to prove that there are infinitely many
$(x,y,z,w)$'s such that the braid index of
$BM_{x,y,z,w}$ is $4$.  More concretely, let ${\cal K}_n :=
BM_{-1, -2, n, 2}$ and we will show that for all $m \geq 1$ the
braid index of ${\cal K}_{2m}$ is $4$. Note that ${\cal K}_2 =
10_{150}$ and ${\cal K}_{2m}$ is a knot.

The {\em enhanced Milnor Number} $\lambda$ defined by Neumann and Rudolph
\cite{NR} is an invariant of fibred knots and links counting the
number of negative Hopf band plumbing to get the fibre surface.
(Recall that the fiber surface of a fibre knot is obtained by
plumbing and deplumbing Hopf bands \cite{Giroux}.)
\begin{lemma}\label{lambda=1}
All ${\cal K}_n$ $(n \geq 2)$ are fibred and have the enhanced
Milnor number  $\lambda = 1.$
\end{lemma}

\begin{figure}[htp]
\begin{center}
\includegraphics [height=150mm]{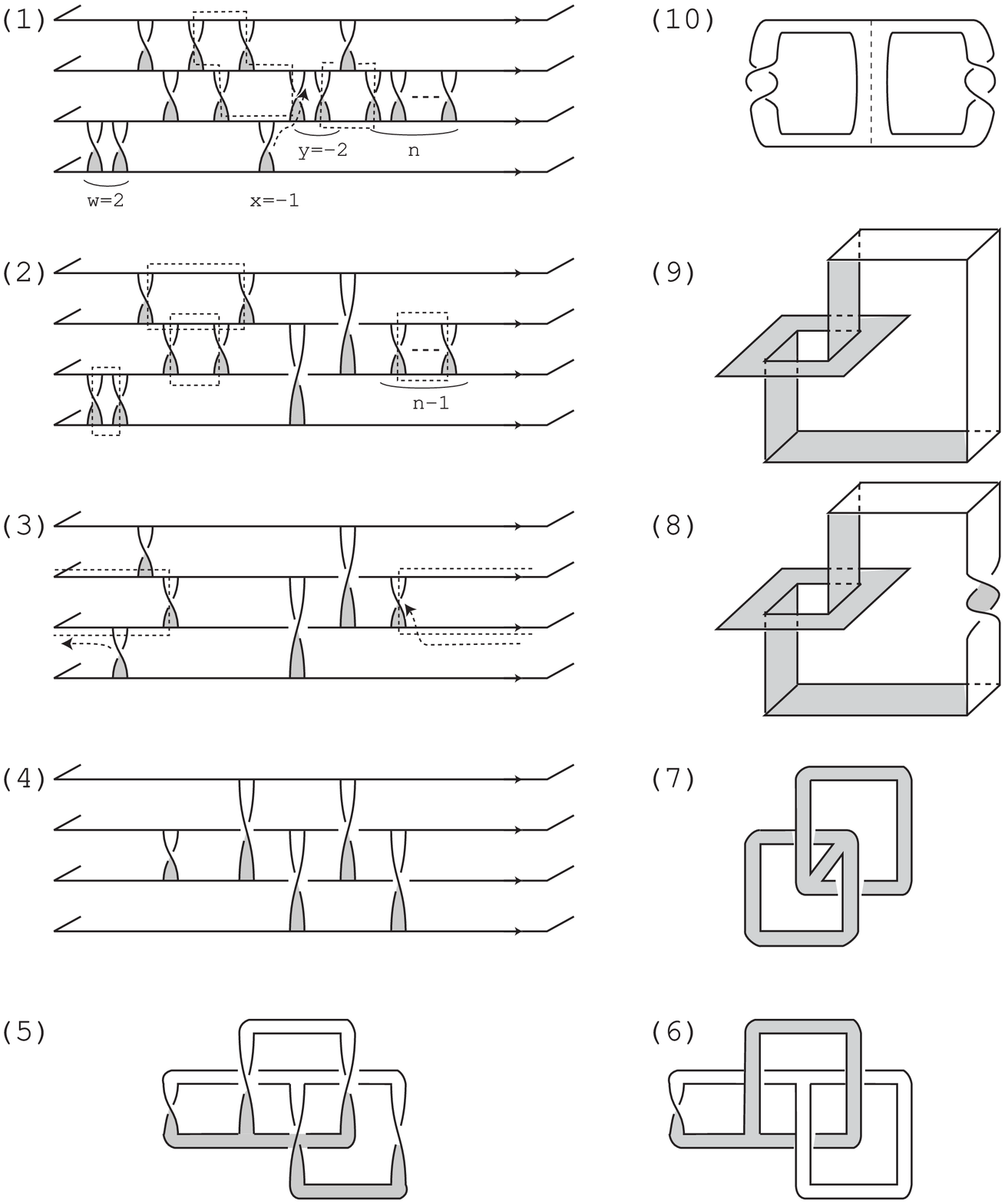}
\end{center}
\caption{} \label{deplum}
\end{figure}

\textbf{Proof of Lemma \ref{lambda=1}. \ } As in the passage $(1)
\Rightarrow (2)$ of Figure \ref{deplum}, we compress twice the
standard Bennequin surface (sketch (1)) of ${\cal K}_n$.
Next, deplumb positive Hopf bands as much as possible as in the
passage sketch $(2) \Rightarrow (3) \Rightarrow (4)=(5)$.  Then
isotope the Seifert surface until we get $P(-2, -2, 2)$ (sketch (8)) a
Pretzel link. These operations do not change the enhanced Milnor
number.

We apply a trick of Melvin and Morton \cite{MM}, as in the passage
sketch $(8) \Rightarrow (9)$ and get $P(-2, 0, 2)$ a Pretzel link.
We remark that the enhanced Milnor number is invariant under this
trick.

Since $P(-2, 0, 2)$ is obtained by plumbing one positive Hopf band
and one negative Hopf band (see sketch (10)), it has the
enhanced Milnor number $\lambda = 1$ so does ${\cal K}_n$.
\hfill $\Box$ %

\begin{figure}[htpb!]
\begin{center}
\includegraphics [height=24mm]{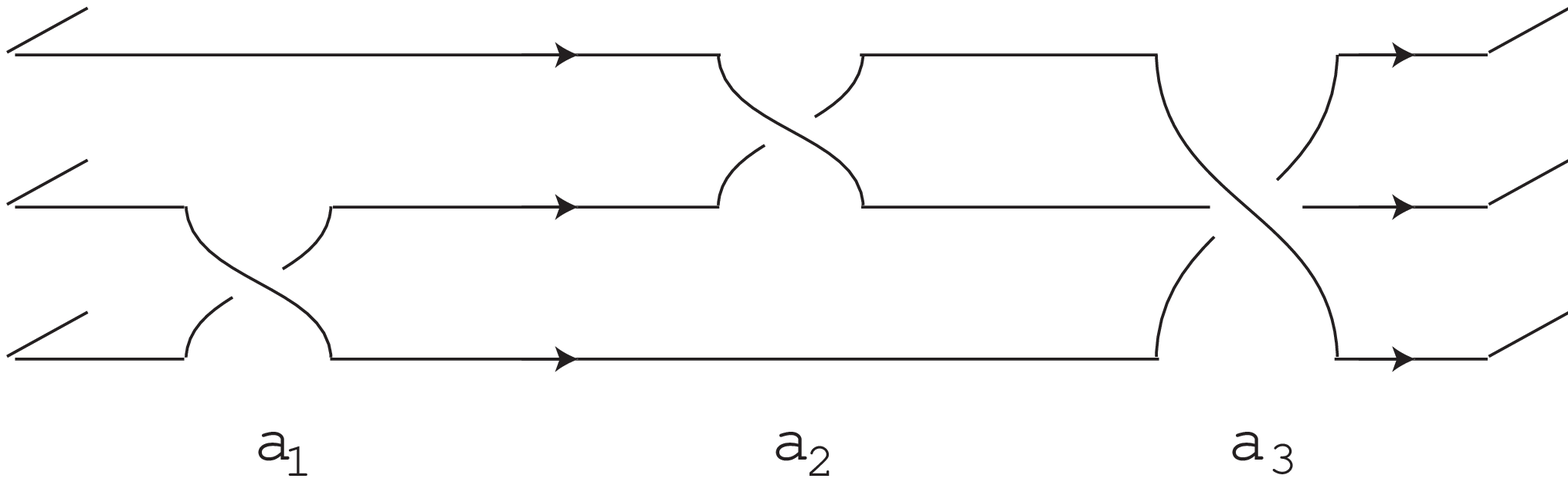}
\end{center}
\caption{} \label{bands}%
\end{figure}

Here we summarize Xu's classification of $3$-braids \cite{Xu}. Let
$\sigma_1, \sigma_2$ be the standard generators of $B_3$ the braid
group of $3$-strings satisfying $\sigma_1 \sigma_2 \sigma_1 =
\sigma_2 \sigma_1 \sigma_2$. Let $a_1 := \sigma_1, a_2 :=
\sigma_2$ and $a_3 := \sigma_2 \sigma_1 \sigma_2^{-1}$. We can
identify them with the twisted bands in Figure \ref{bands}. Let
$\alpha := a_1 a_3 = a_2 a_1 = a_3 a_2$. If $w \in B_3$ let
$\overline{w}$ denote $w^{-1}.$

\begin{theorem}{\rm (Xu \cite{Xu}.)}
Every conjugacy class in $B_3$ can be represented by a shortest
word in $a_1, a_2, a_3$ uniquely up to symmetry. And the word has
one of the three forms:
$$(1) \alpha^k P, \quad (2) N \overline{\alpha}^{k} , \quad (3) NP. $$
where $k\geq 0$ and $\overline{N}, P$ are positive words and the
arrays of subscripts of the words are non-decreasing.
\end{theorem}

%Thanks to Xu's result we have:
\begin{lemma}\label{ABCD}
If a closed $3$-braid has $\lambda=1$ and is a knot, then up to
symmetry it has one of the following Xu's forms:
\begin{eqnarray*}
A_x &:=& \overline{a_3} \ \overline{a_2} \ (a_1)^x, \quad x\geq 2, \mbox{\ even,} \\%
B_{x, y} &:=& \overline{a_3} \ \overline{a_3} \ (a_1)^x (a_2)^y,
                                     \quad x, y \geq 3, \mbox{\ odd,} \\%
C_{x, y, z} &:=& \overline{a_2} \ (a_1)^x (a_2)^y (a_3)^z,
                                     \quad x+z= \mbox{odd, } \ y= \mbox{even, } \ x,y,z\geq 1,\\ %
D_{x, y, z, w} &:=& \overline{a_2} \ (a_1)^x (a_2)^y (a_3)^z (a_1)^w, %
                         \quad x,y \geq 2, \ z,w\geq 1. %
\end{eqnarray*}
\end{lemma}

\textbf{Proof of Lemma \ref{ABCD}.} For simplicity
let $\longrightarrow$ (resp. $\Longrightarrow$) denote
``deplumbing of positive-Hopf (resp. negative) bands''.
We denote $w=w'$ when $w, w'$ have the same conjugacy class.
Assume we have a word $w \in B_3.$

\begin{figure}[htpb!]
\begin{center}
\includegraphics [height=24mm]{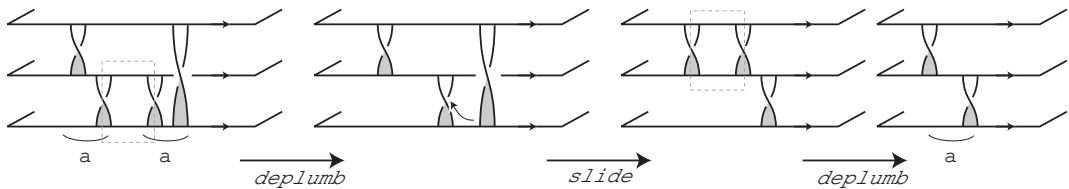}
\end{center}
\caption{$\alpha^2 \longrightarrow \alpha.$} \label{alpha}%
\end{figure}

\textbf{Case (1)-1}.
Suppose $w=\alpha^k$ for some $k\geq 1.$
Since $\alpha^2 \longrightarrow \alpha$ (see Figure \ref{alpha}),
we have $w=\alpha^k \longrightarrow \alpha.$ Since the braid closure of
$\alpha$ is the unknot, $w$ has $\lambda=0.$

\begin{figure}[htpb!]
\begin{center}
\includegraphics [height=24mm]{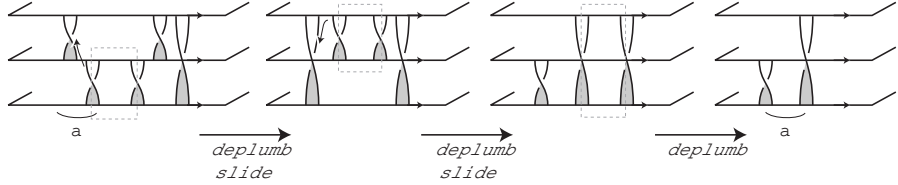}
\end{center}
\caption{$\alpha a_1 a_2 a_3 \longrightarrow \alpha.$}
\label{alpha2}\end{figure}

\textbf{Case (1)-2}. If $w=\alpha^k P$ ($k\geq 1$), up to
permutation of subscripts $\{1,2,3\}$ we get
$$\alpha^k P \longrightarrow \alpha P \longrightarrow \alpha (a_1 a_2
a_3 a_1 a_2 a_3 \cdots \cdots).$$
Since $\alpha a_1 a_2 a_3 \longrightarrow \alpha$ (see
Figure \ref{alpha2}) we have $$\alpha \ \overbrace{a_1 a_2
a_3 a_1 a_2 a_3 \cdots \cdots}^{\mbox{length$=l$}} \
\longrightarrow \ \alpha \ \overbrace{a_1 a_2 a_3 a_1 a_2 a_3
\cdots \cdots}^{\mbox{length$=l-3$}} \quad \mbox{for } l\geq 3.$$
If $l=1, 2$, we have $\alpha a_1 \longrightarrow \alpha$ and $\alpha a_1 a_2
\longrightarrow \alpha$. Thus $w$ has $\lambda=0.$

\textbf{Case (1)-3}. Assume $w=P$ (with
no $\alpha$ part). There are three possible cases to study:
$$P \longrightarrow (a_1 a_2 a_3)^n, \quad
P \longrightarrow (a_1 a_2 a_3)^n a_1 \ \mbox{ and } \
P \longrightarrow (a_1 a_2 a_3)^n a_1 a_2.$$

If $P$ satisfies the third case, since $(a_1 a_2 a_3)^n a_1 a_2
=  a_2 (a_1 a_2 a_3)^n a_1 = \alpha (a_2 a_3 a_1)^n \longrightarrow \alpha,$
$w$ has $\lambda=0.$

If $P$ satisfies the second case, since $(a_1 a_2 a_3)^n a_1 =
a_1 (a_1 a_2 a_3)^n \longrightarrow (a_1 a_2 a_3)^n$ this case can be
reduced to the first case.

If $P$ satisfies the first case, it is known that closure of $P$ is not fibred
\cite{H}\cite{S} i.e., $w$ is not fibred.

\textbf{Case (2)-1}.
Assume $w=\overline{\alpha}^k$ for some $k\geq 1.$ Figure \ref{alpha} shows
that $\overline{\alpha}^2
\Longrightarrow \overline{\alpha}$ by deplumbing negative-Hopf band
twice, i.e., $\overline{\alpha}^2$ has $\lambda=2$. Thus
$\overline{\alpha}^k$ has $\lambda=2(k-1)\neq 1$.

\textbf{Case (2)-2}.
Suppose $w=N \overline{\alpha}^k$ where $k\geq 1$.

If $w=\overline{a_i}\ \overline{\alpha}$ we have
$\overline{a_i}\ \overline{\alpha} \Longrightarrow \overline{\alpha}$ and
$w$ has $\lambda=1.$ However, the closure of $w$ has more than one component
and it does not satisfy the condition of the lemma.

If $w\neq \overline{a_i}\ \overline{\alpha}$
by similar argument as in case (1)-2, we have $N
\overline{\alpha}^k \Longrightarrow \overline{\alpha}$ and $w$
has $\lambda \geq 2.$

\textbf{Case (2)-3}.
Suppose $w=N$ (no $\overline{\alpha}$ part).

Assume $N \Longrightarrow (\overline{a_3} \ \overline{a_2}\ \overline{a_1})^n
\overline{a_3} \ \overline{a_2}$ for some $n\geq 0.$
If $n=0$ then $w$ has $\lambda=1$ if and only if
$w=\overline{a_3}\ \overline{a_3} \ \overline{a_2}.$ However it has two
components and does not satisfy the condition of the lemma.
If $n\geq 1,$ since
$(\overline{a_3} \ \overline{a_2} \ \overline{a_1})^n
  \overline{a_3} \ \overline{a_2}
= \overline{a_2} \ \overline{a_3}
  (\overline{a_2} \ \overline{a_1} \ \overline{a_3})^n
=\overline{\alpha}\ (\overline{a_2} \ \overline{a_1} \ \overline{a_3})^n
\Longrightarrow \overline{\alpha},$ $w$ has $\lambda \geq 3n.$

Assume $N \Longrightarrow (\overline{a_3} \ \overline{a_2}\ \overline{a_1})^n
\overline{a_3}$ for some $n\geq 0.$
If $n=0$ then  $w$ has $\lambda=1$ if and only if
$w=\overline{a_3}\ \overline{a_3}$. However this has two components.
If $n\geq 1,$ since $(\overline{a_3} \ \overline{a_2}\ \overline{a_1})^n
\overline{a_3} \Longrightarrow
(\overline{a_3} \ \overline{a_2}\ \overline{a_1})^n$ it can be reduced to the
next case we discuss.

Assume $N \Longrightarrow (\overline{a_3} \ \overline{a_2}\
\overline{a_1})^n$ then it is known that $w$ is not fibred \cite{H}\cite{S}.

\textbf{Case (3)}. Assume $w=NP$ for some $N, P \neq \emptyset.$

We introduce new symbol ``$\approx$'' denoting Melvin and Morton's
trick \cite{MM}. In our situation we have
$$
\overline{a_i} a_{i-1} a_i \approx \overline{a_i}\
\overline{a_{i-1}} a_i\ \mbox{ and }\ a_i a_{i+1} \overline{a_i} \approx
a_i \overline{a_{i+1}}\ \overline{a_i}.$$
Recall that this
trick does {\em not} change $\lambda$ nor fibre-ness.

Let $\leadsto$ denote composition of $\pm$ Hopf bands deplumbings.

After deplumbing $\pm$ Hopf bands sufficiently enough times, $w$ can be
reduced to one
of the following 18 forms up to permutation of $\{1,2,3\}$.

\begin{tabular}{|l|ll|}
\hline case & word NP &  \\ \hline
i & $(\overline{a_{2}}\ \overline{a_{1}}\ \overline{a_{3}}%
)^{k}(a_{1}a_{2}a_{3})^{l}$ & $k\geq 1,l\geq 1$ \\ \hline
ii & $(\overline{a_{2}}\ \overline{a_{1}}\ \overline{a_{3}}%
)^{k}(a_{1}a_{2}a_{3})^{l}a_{1}$ & $k\geq 1,l\geq 0$ \\ \hline
iii & $(\overline{a_{2}}\ \overline{a_{1}}\ \overline{a_{3}}%
)^{k}(a_{1}a_{2}a_{3})^{l}a_{1}a_{2}$ & $k\geq 1,l\geq 0$ \\
\hline
iv & $\overline{a_{3}}(\overline{a_{2}}\ \overline{a_{1}}\ \overline{a_{3}}%
)^{k}(a_{1}a_{2}a_{3})^{l}$ & $k\geq 0,l\geq 1$ \\ \hline
v & $\overline{a_{3}}(\overline{a_{2}}\ \overline{a_{1}}\ \overline{a_{3}}%
)^{k}(a_{1}a_{2}a_{3})^{l}a_{1}$ & $k\geq 0,l\geq 0$ \\ \hline
vi & $\overline{a_{3}}(\overline{a_{2}}\ \overline{a_{1}}\ \overline{a_{3}}%
)^{k}(a_{1}a_{2}a_{3})^{l}a_{1}a_{2}$ & $k\geq 0,l\geq 0$ \\
\hline
vii & $\overline{a_{1}}\ \overline{a_{3}}(\overline{a_{2}}\ \overline{a_{1}}%
\ \overline{a_{3}})^{k}(a_{1}a_{2}a_{3})^{l}$ & $k\geq 0,l\geq 1$
\\ \hline
viii & $\overline{a_{1}}\ \overline{a_{3}}(\overline{a_{2}}\ \overline{a_{1}}%
\ \overline{a_{3}})^{k}(a_{1}a_{2}a_{3})^{l}a_{1}$ & $k\geq 0,l\geq 0$ \\
\hline ix & $\overline{a_{1}}\ \overline{a_{3}}(\overline{a_{2}}\
\overline{a_{1}}\
\overline{a_{3}})^{k}(a_{1}a_{2}a_{3})^{l}a_{1}a_{2}$ & $k\geq 0,l\geq 0$ \\
\hline
\end{tabular}%
\begin{tabular}{|l|ll|}
\hline case & word NP &  \\ \hline
i' & $(\overline{a_{1}}\ \overline{a_{3}}\ \overline{a_{2}}%
)^{k}(a_{1}a_{2}a_{3})^{l}$ & $k\geq 1,l\geq 1$ \\ \hline
ii' & $(\overline{a_{1}}\ \overline{a_{3}}\ \overline{a_{2}}%
)^{k}(a_{1}a_{2}a_{3})^{l}a_{1}$ & $k\geq 1,l\geq 0$ \\ \hline
iii' & $(\overline{a_{1}}\ \overline{a_{3}}\ \overline{a_{2}}%
)^{k}(a_{1}a_{2}a_{3})^{l}a_{1}a_{2}$ & $k\geq 1,l\geq 0$ \\
\hline
iv' & $\overline{a_{2}}(\overline{a_{1}}\ \overline{a_{3}}\ \overline{a_{2}}%
)^{k}(a_{1}a_{2}a_{3})^{l}$ & $k\geq 0,l\geq 1$ \\ \hline
v' & $\overline{a_{2}}(\overline{a_{1}}\ \overline{a_{3}}\ \overline{a_{2}}%
)^{k}(a_{1}a_{2}a_{3})^{l}a_{1}$ & $k\geq 0,l\geq 0$ \\ \hline
vi' & $\overline{a_{2}}(\overline{a_{1}}\ \overline{a_{3}}\ \overline{a_{2}}%
)^{k}(a_{1}a_{2}a_{3})^{l}a_{1}a_{2}$ & $k\geq 0,l\geq 0$ \\
\hline
vii' & $\overline{a_{3}}\ \overline{a_{2}}(\overline{a_{1}}\ \overline{a_{3}}%
\ \overline{a_{2}})^{k}(a_{1}a_{2}a_{3})^{l}$ & $k\geq 0,l\geq 1$
\\ \hline
viii' & $\overline{a_{3}}\ \overline{a_{2}}(\overline{a_{1}}\ \overline{a_{3}%
}\ \overline{a_{2}})^{k}(a_{1}a_{2}a_{3})^{l}a_{1}$ & $k\geq 0,l\geq 0$ \\
\hline
ix' & $\overline{a_{3}}\ \overline{a_{2}}(\overline{a_{1}}\ \overline{a_{3}}%
\ \overline{a_{2}})^{k}(a_{1}a_{2}a_{3})^{l}a_{1}a_{2}$ & $k\geq
0,l\geq 0$
\\ \hline
\end{tabular}

For example, assume $w$ can be reduced to have form iv'.

Assume $k=0, l\geq 1$ i.e.,
$w \longrightarrow \overline{a_{2}} (a_{1}a_{2}a_{3})^l.$ Since
$$\overline{a_{2}} a_{1}a_{2}a_{3} \approx \overline{a_{2}}\ \overline{a_{1}}
a_{2}a_{3}  = \overline{a_{2}} a_{3}a_{3}\overline{a_{2}}
\longrightarrow \overline{a_{2}} a_{3}\overline{a_{2}} = a_{1}
\overline{a_{2}} \ \overline{a_{2}} \Longrightarrow a_{1}
\overline{a_{2}}  = \mbox{ unknot,}$$
$\overline{a_{2}} (a_{1}a_{2}a_{3})^l$ has $\lambda = 1$ if and
only if $l=1.$
Let $$C_{x,y,z}:=\overline{a_{2}} (a_{1})^x (a_{2})^y (a_{3})^z
\mbox{ for $x,y,z\geq 1.$}$$
Since $C_{x,y,z} \longrightarrow \overline{a_{2}} a_{1}a_{2}a_{3}$,
$C_{x,y,z}$ has $\lambda=1.$

To study rest of the cases ($k,l\geq 1$) we remark that
$(\overline{a_{1}}\ \overline{a_{3}}\
\overline{a_{2}})^{k}(a_{1}a_{2}a_{3})^{k}$ can be reduced to $\overline{a_1}
a_3$ by deplumbing positive and negative Hopf bands
each $(3k-1)$-times i.e., $(\overline{a_{1}}\ \overline{a_{3}}\
\overline{a_{2}})^{k}(a_{1}a_{2}a_{3})^{k} \leadsto \overline{a_1}
a_3.$

If $k=l \geq 1$ then
$$w \leadsto
\overline{a_{2}}(\overline{a_{1}}\ \overline{a_{3}}\ \overline{a_{2}}%
)^{k}(a_{1}a_{2}a_{3})^{k} \ \leadsto \ \overline{a_{2}}
(\overline{a_1} a_3) = a_1 \overline{a_2}\ \overline{a_2}
\Longrightarrow a_1 \overline{a_2} = \mbox{ unknot }$$ and $w$ has
$\lambda \geq 2.$

If $k>l \geq 1$ then
\begin{eqnarray*}
 \lefteqn{w \leadsto \overline{a_{2}}(\overline{a_{1}}\ \overline{a_{3}}\
  \overline{a_{2}} )^{k}(a_{1}a_{2}a_{3})^{l} \leadsto
  \overline{a_{2}}(\overline{a_{1}}\ \overline{a_{3}}\
  \overline{a_{2}})^{k-l}(\overline{a_1} a_3)  = (\overline{a_{1}}\
  \overline{a_{3}}\ \overline{a_{2}})^{k-l} \overline{a_1} a_3
  \overline{a_{2}} } \\ %
 && \Longrightarrow (\overline{a_{1}}\ \overline{a_{3}}\
  \overline{a_{2}})^{k-l} \overline{a_1} a_3 = (\overline{a_{1}}\
  \overline{a_{3}}\ \overline{a_{2}})^{k-l-1} \overline{a_1}\
  \overline{a_3} a_1 \overline{a_{2}}\ \overline{a_{2}}
  \Longrightarrow \approx (\overline{a_{1}}\ \overline{a_{3}}\
  \overline{a_{2}})^{k-l-1} \overline{a_1} a_3 a_1 \overline{a_{2}} \\ %
 && =  (\overline{a_{1}}\ \overline{a_{3}}\ \overline{a_{2}})^{k-l-1}
  \overline{a_1}\ \overline{a_1}\ a_3 a_3 \leadsto
  (\overline{a_{1}}\ \overline{a_{3}}\
  \overline{a_{2}})^{k-l-1} \overline{a_1} a_3 \leadsto \overline{a_1} a_3
  = \mbox{ unknot }
\end{eqnarray*}
and $w$ has $\lambda\geq 2.$

If $l>k \geq 1$ then
\begin{eqnarray*}
 w &\leadsto & x\overline{a_{2}}(\overline{a_{1}}\ \overline{a_{3}}\
  \overline{a_{2}} )^{k}(a_{1}a_{2}a_{3})^{l} \leadsto
  \overline{a_{2}}(\overline{a_1} a_3)(\overline{a_{1}}\ \overline{a_{3}}\
  \overline{a_{2}})^{l-k} \Longrightarrow a_1  \overline{a_{2}}
  (\overline{a_{1}}\ \overline{a_{3}}\ \overline{a_{2}})^{l-k} \\ %
 &\approx \leadsto & a_1  \overline{a_{2}}
  (\overline{a_{1}}\ \overline{a_{3}}\ \overline{a_{2}})^{l-k-1}
  \approx \leadsto a_1  \overline{a_{2}} = \mbox{ unknot } %
\end{eqnarray*}
and $w$ has $\lambda\geq 2.$

Thus if $w=NP$ for some $N,P$ and can be reduced to have form iv'
then $w$ has $\lambda=1$ if and only if $w=C_{x,y,z}$ for
$x,y,z\geq 1$. To make the braid closure of $w$ have one component, we
further require $x+z= \mbox{odd.}$

The following table lists all the words with $\lambda=1.$

\begin{tabular}{|l|l|}
\hline case & word with $\lambda =1.$ \\ \hline %
i & none. \\ \hline %
ii & $\overline{a_{2}}\ \overline{a_{1}}\
\overline{a_{3}}\ a_{1}^{x}\ $(2 or 3 components.) \\ \hline %
iii & reduced to Case (1) or (2). \\ \hline %
iv & reduced to iii. \\ \hline %
v &
\begin{tabular}{l}
$\overline{a_{3}}\ a_{1}^{x}\ a_{2}^{y}\ a_{3}^{z}\ a_{1}^{w}\ =\{%
\begin{tabular}{l}
$C_{x+1,y,z}\ $Case iv' when $w=1,$ \\
$D_{x+1,y,z,w-1}$\ Case v' when $w\geq 2.$%
\end{tabular}%
$ \\
$\overline{a_{3}}\ \overline{a_{3}}\ a_{1}^{x}\ $(2 or 3 components.)%
\end{tabular}
\\ \hline%
vi &
\begin{tabular}{l}
$\overline{a_{3}}\ \overline{a_{3}}\ a_{1}^{x}\ a_{2}^{y}=:B_{x,y}.$ \\ %
$\overline{a_{3}}\ a_{1}^{x}\ a_{2}^{y}\ a_{3}^{z}\ a_{1}^{w}\ a_{2}^{v}\ =\{%
\begin{tabular}{l}
$C_{x+v+1,y,z}\ $Case iv' when $w=1,$ \\
$D_{x+v+1,y,z,w-1}$\ Case v' when $w\geq 2.$%
\end{tabular}%
$%
\end{tabular}
\\ \hline %
vii & $\overline{a_{1}}\ \overline{a_{3}}\ a_{1}^{x}\ a_{2}^{y}\ a_{3}^{z}\ =%
\overline{a_{1}}\ \overline{a_{3}}\ a_{1}^{x+z}\ a_{2}^{y}\ $Case ix. \\ \hline %
viii & reduced to iv. \\ \hline %
ix & $\overline{a_{1}}\ \overline{a_{3}}\ a_{1}^{x}\ a_{2}^{y}\
=B_{x+1,y-1}$ Case v or vi. \\ \hline i' & none. \\ \hline ii' &
reduced to Case (1) or (2). \\ \hline iii' & none. \\ \hline
iv' & $\overline{a_{2}}\ a_{1}^{x}\ a_{2}^{y}\ a_{3}^{z}\ =:C_{x,y,z}.$ \\
\hline v' &
\begin{tabular}{l}
$\overline{a_{2}}\ a_{1}^{x}\ a_{2}^{y}\ a_{3}^{z}\ a_{1}^{w}\
=:D_{x,y,z,w}. $ \\
$\overline{a_{2}}\ \overline{a_{2}}\ a_{1}^{x}\ $(2 or 3 components.)%
\end{tabular}
\\ \hline
vi' & reduced to ii' \\ \hline vii' & reduced to Case (1) or (2).
\\ \hline viii' & $\overline{a_{3}}\ \overline{a_{2}}\ a_{1}^{x}.\
=:A_{x}.$ \\ \hline
ix' & $\overline{a_{3}}\ \overline{a_{2}}\ a_{1}^{x}\ a_{2}^{y}$ $=\{%
\begin{tabular}{l}
$\overline{a_{3}}\ \overline{a_{3}}\ a_{2}^{y+1}\ $Case v' when $x=1,$ \\
$B_{x-1,y+1}$\ Case v' when $x\geq 2.$%
\end{tabular}%
$ \\ \hline
\end{tabular}

Words $A_x, \cdots, D_{x,y,z,w}$ are defined as above.  Table shows that any
$w$ with $\lambda=1$ and having one component has one of the forms;
$A_x, \cdots, D_{x,y,z,w}.$

\hfill $\Box$ %

\begin{lemma}\label{alexander}
Leading terms of Alexander polynomials of ${\cal K}_{n}$, $A_x$,
$B_{x, y}$, $C_{x, y, z}$ and $D_{x, y, z, w}$ are the following:
\begin{eqnarray*}
{\cal K}_n; && \pm(1 - 4t -6t^2 + 8t^3 - \cdots) \quad \mbox{if $n \geq 2$,} \\
A_x;        && \pm(1 - 3t + \cdots) \quad \mbox{if } x\geq 2, \\%
B_{x, y};   && \pm(1 - 3t + \cdots) \quad \mbox{if } x, y \geq 3, \\%
C_{x, y, z}; && \pm(1 - 5t + \cdots) \quad \mbox{if $x,z \geq 2$,}\\ %
C_{1, 2, z},  C_{1, y, 2}, C_{2, y, 1}, C_{x, 2, 1}; %
            && \pm(1 - 4t + 6t^2 -7t^3 + \cdots) %
            \quad \mbox{if $x, y, z \geq 4$,}\\ %
C_{1, y, z}, C_{x, y, 1}; %
            && \pm(1 - 4t + 7t^2 + \cdots) \quad \mbox{if $x, y, z \geq 3$}\\%
D_{x, y, z, w}, D_{x, y, z, 1}; && \pm(1 - 6t + \cdots) \quad \mbox{if } x,y,z,w \geq 2, \\
D_{x, y, 1, w}; && \pm(1 - 5t + \cdots) \quad \mbox{if } x,y,w \geq 2. %
\end{eqnarray*}
\end{lemma}

\begin{figure}[htpb!]
\begin{center}
\includegraphics [height=45mm]{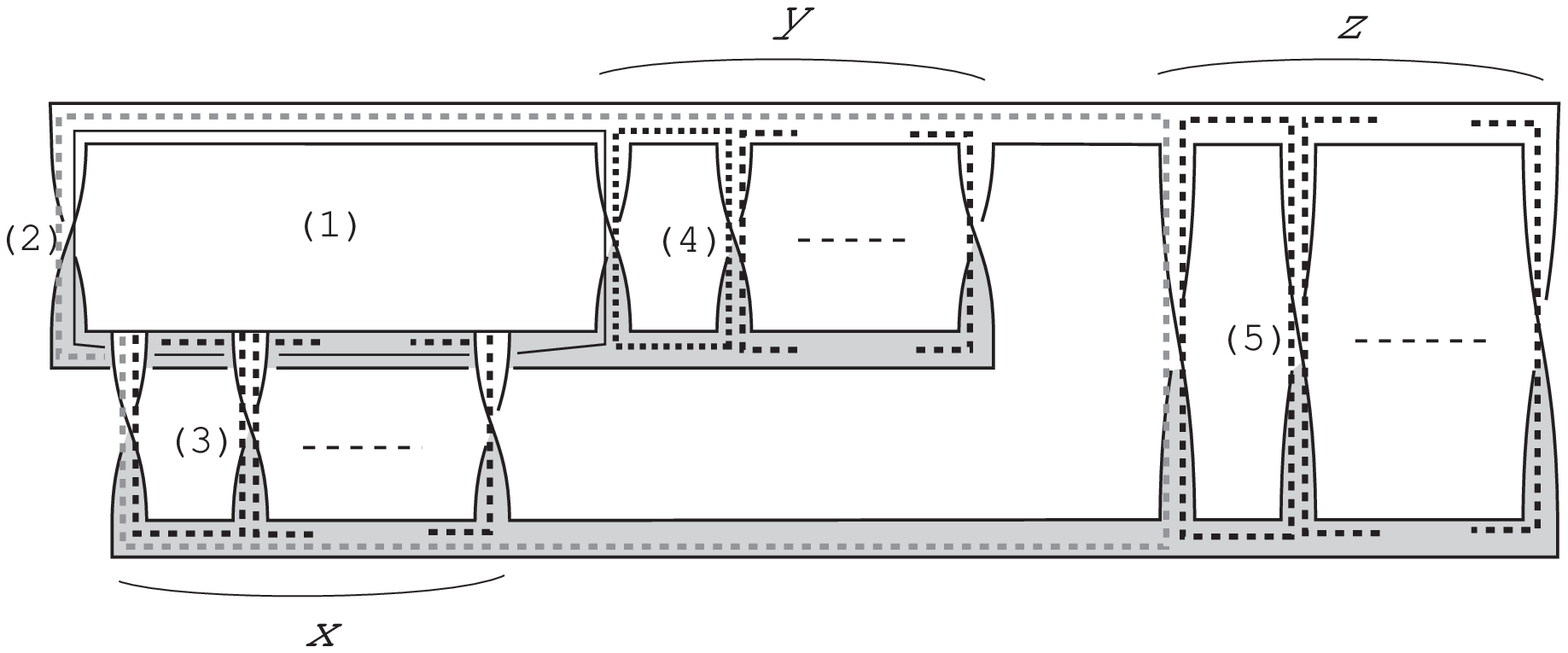}
\end{center}
\caption{The Bennequin surface $F$ of $C_{x, y, z}=\overline{a_2} \ (a_1)^x (a_2)^y (a_3)^z$
         and a basis for $H_1(F).$} \label{2-123}%
\end{figure}

\textbf{Proof of Lemma \ref{alexander}.} We prove that the
Alexander polynomial of $C_{x, y, z}$ for some $x, y, z \geq 2$ is
$\pm(1 - 5t + \cdots)$.  Recall that the Bennequin surface of Xu's
form gives a minimal genus Seifert surface. Let $F$ be the
Bennequin surface of $C_{x, y, z}$  and choose a basis
$$\{ u^{(1)}, u^{(2)}, u^{(3)}_1, \cdots, u^{(3)}_{x-1}, u^{(4)}_1,
\cdots, u^{(4)}_{y-1}, u^{(5)}_1, \cdots, u^{(5)}_{z-1} \}$$ for
$H_1(F)$ as in Figure \ref{2-123}. In the sketch, $u^{(k)}$
($k=1,2,3,4,5$) corresponds to the loop $(k)$.
With respect to the basis, let $V_{x,y,z}$ denote the Seifert matrix for
$C_{x, y, z}$.
$$ V_{x,y,z} = \left[
  \begin{array}{c|c|cccc|cccc|cccc}
& 1 &  &  &  &  & 1 &  &  &  &  &  &  &  \\ \hline %
1 &  & -1 &  & &  &  &  &  &  & 1 &  &  &  \\ \hline%
&  & -1 & 1 &  &  &  &  &  &  &  &  &  &  \\
&  &  & -1 & \ddots  &  &  &  &  &  &  &  &  &  \\
&  &  &  & \ddots  & 1 &  &  &  &  &  &  &  &  \\
&  &  &  &  & -1 &  &  &  &  &  &  &  &  \\ \hline%
&  &  &  &  &  & -1 & 1 &  &  &  &  &  &  \\
&  &  &  &  &  &  & -1 & \ddots  &  &  &  &  &  \\
&  &  &  &  &  &  &  & \ddots  & 1 &  &  &  &  \\
&  &  &  &  &  &  &  &  & -1 &  &  &  &  \\ \hline%
&  &  &  &  &  &  &  &  &  & -1 & 1 &  &  \\
&  &  &  &  &  &  &  &  &  &  & -1 & \ddots  &  \\
&  &  &  &  &  &  &  &  &  &  &  & \ddots  & 1 \\
&  &  &  &  &  &  &  &  &  &  &  &  & -1%
\end{array} \right] $$
It has $0$'s in all the blank places. The $3$rd (resp. $4$th,
$5$th) diagonal block has size $(x-1)\times (x-1)$ (resp.
$(y-1)\times (y-1)$, $(z-1)\times (z-1)$). Alexander polynomial
satisfies:

$\Delta_{x,y,z}(t) = \det(V_{x,y,z}^T - tV_{x,y,z})$
$$
= \det \left[
\begin{array}{c|c|cccc|cccc|cccc}
& 1-t &  &  &  &  & -t &  &  &  &  &  &  &  \\ \hline%
1-t &  & t &  &  &  &  &  &  &  & -t &  &  &  \\ \hline%
& -1 & -1+t & -t &  &  &  &  &  &  &  &  &  &  \\
&  & 1 & \ddots  & \ddots  &  &  &  &  &  &  &  &  &  \\
&  &  & \ddots  & \ddots  & -t &  &  &  &  &  &  &  &  \\
&  &  &  & 1 & -1+t &  &  &  &  &  &  &  &  \\ \hline%
1 &  &  &  &  &  & -1+t & -t &  &  &  &  &  &  \\
&  &  &  &  &  & 1 & \ddots  & \ddots  &  &  &  &  &  \\
&  &  &  &  &  &  & \ddots  & \ddots  & -t &  &  &  &  \\
&  &  &  &  &  &  &  & 1 & -1+t &  &  &  &  \\ \hline%
& 1 &  &  &  &  &  &  &  &  & -1+t & -t &  &  \\
&  &  &  &  &  &  &  &  &  & 1 & \ddots  & \ddots  &  \\
&  &  &  &  &  &  &  &  &  &  & \ddots  & \ddots  & -t \\
&  &  &  &  &  &  &  &  &  &  &  & 1 & -1+t%
\end{array} \right]. $$
Expanding it in the $(x+1)$th column, we have

$ \Delta_{x,y,z}(t) = (-1+t)\Delta_{x-1,y,z}(t)$
$$ -(-t)\det
\left[
\begin{array}{c|c|cccc|ccc|ccc}
& 1-t &  &  &  &  & -t &  & &  &  &  \\ \hline%
1-t &  & t &  &  &  &  &  &  & -t &  &  \\ \hline%
& -1 & -1+t & -t &  &  &  &  &  &  &  &  \\
&  & 1 & \ddots  & -t &  &  &  &  &  &  &  \\
&  &  & 1 & -1+t & -t &  &  &  &  &  &  \\
&  &  &  &  & 1 &  &  &  &  &  &  \\ \hline%
1 &  &  &  &  &  & -1+t & -t &  &  &  &  \\
&  &  &  &  &  & 1 & \ddots  & -t &  &  &  \\
&  &  &  &  &  &  & 1 & -1+t &  &  &  \\ \hline%
& 1 &  &  &  &  &  &  &  & -1+t & -t &  \\
&  &  &  &  &  &  &  &  & 1 & \ddots  & -t \\
&  &  &  &  &  &  &  &  &  & 1 & -1+t%
\end{array} \right] $$
$= (-1+t)\Delta_{x-1,y,z}(t) + t \Delta_{x-2,y,z}(t).$

If $\Delta_{i,y,z}(t)=(-1)^i (\alpha_0 + \alpha_1 t + \alpha_2 t^2
+ \cdots )$ for $i=x-1$ and $x-2$, then
\begin{eqnarray*}
\Delta_{x,y,z}(t) &=& (-1+t) (-1)^{x-1} (\alpha_0 + \alpha_1 t +
\alpha_2 t^2 + \cdots )
+ t (-1)^{x-2} (\alpha_0 + \alpha_1 t + \alpha_2 t^2 + \cdots) \\ %
&=& (-1)^x (\alpha_0 + \alpha_1 t + \alpha_2 t^2 + \cdots ).
\end{eqnarray*}
In fact, $\Delta_{x,y,z}(t) = (-1)^{x+y+z} (1-5t+ \cdots )$ for
all $x,y,z \in \{2,3\}$.  By induction, $\Delta_{x,y,z}(t) =
(-1)^{x+y+z} (1-5t+ \cdots )$ for all $x, y, z \geq 2.$

Other cases follow by similar arguments.
\hfill $\Box$ %

\textbf{Proof of Theorem \ref{BM-thm}.} By Lemmas
\ref{lambda=1}, \ref{ABCD}, \ref{alexander}, our knot ${\cal
K}_{2m}$ where ($m\geq 1$) cannot be a $3$-braid. Then by Lemma
\ref{D+}, Theorem \ref{BM-thm} follows.
\hfill $\Box$ %

Department of Mathematics, Columbia University, New York, NY 10027 \\
{\it E-mail address:} {\tt yuri@math.columbia.edu}
\end{document}